\documentclass[reqno]{amsart}

\usepackage[latin1]{inputenc}
\usepackage{amssymb}
\usepackage{graphicx}
\usepackage{amscd}
\usepackage[hidelinks]{hyperref}
\usepackage{color}
\usepackage{float}
\usepackage{graphics,amsmath,amssymb}
\usepackage{amsthm}
\usepackage{amsfonts}
\usepackage{latexsym}
\usepackage{epsf}
\usepackage{enumerate}
\usepackage{xifthen}
\usepackage{mathrsfs}
\usepackage{dsfont}
\usepackage{makecell}
\usepackage[FIGTOPCAP]{subfigure}
\usepackage{amsmath}
\allowdisplaybreaks[4]
\usepackage{listings}
\usepackage{etoolbox}
\usepackage{fancyhdr}

 \newtheoremstyle{mytheorem}
 {3pt}
 {3pt}
 {\slshape}
 {}
 {\bfseries}
 {.}
 { }
 {}

\numberwithin{equation}{section}

\theoremstyle{theorem}
\newtheorem{theorem}{Theorem}[section]

\newtheorem{conjecture}[theorem]{Conjecture}

\theoremstyle{definition}

\newcommand{\doi}[1]{\href{http://dx.doi.org/#1}{doi: #1}}
\renewcommand{\MR}[1]{\href{http://www.ams.org/mathscinet-getitem?mr=#1}{MR#1}}

\newcommand{\Keywords}[1]{\ifthenelse{\isempty{#1}}{}{\smallskip \smallskip \noindent \textbf{Keywords}. #1}}
\newcommand{\MSC}[2][2010]{\ifthenelse{\isempty{#2}}{}{\smallskip \smallskip \noindent \textbf{#1MSC}. #2}}
\newcommand{\abstractnote}[1]{\ifthenelse{\isempty{#1}}{}{\smallskip \smallskip \noindent \textsuperscript{\dag}#1}}

\makeatletter
\def\specialsection{\@startsection{section}{1}%
  \z@{\linespacing\@plus\linespacing}{.5\linespacing}%
  {\normalfont}}
\def\section{\@startsection{section}{1}%
  \z@{.7\linespacing\@plus\linespacing}{.5\linespacing}%
  {\normalfont\scshape}}
\patchcmd{\@settitle}{\uppercasenonmath\@title}{\Large\boldmath}{}{}
\patchcmd{\@settitle}{\begin{center}}{\begin{flushleft}}{}{}
\patchcmd{\@settitle}{\end{center}}{\end{flushleft}}{}{}
\patchcmd{\@setauthors}{\MakeUppercase}{\normalsize}{}{}
\patchcmd{\@setauthors}{\centering}{\raggedright}{}{}
\patchcmd{\section}{\scshape}{\large\bfseries\boldmath}{}{}
\patchcmd{\subsection}{\bfseries}{\bfseries\boldmath}{}{}
\renewcommand{\@secnumfont}{\bfseries}
\patchcmd{\@startsection}{\@afterindenttrue}{\@afterindentfalse}{}{}
\patchcmd{\abstract}{\leftmargin3pc}{\leftmargin1pc}{}{}

\def\maketitle{\par
  \@topnum\z@ 
  \@setcopyright
  \thispagestyle{empty}
  \ifx\@empty\shortauthors \let\shortauthors\shorttitle
  \else \andify\shortauthors
  \fi
  \@maketitle@hook
  \begingroup
  \@maketitle
  \toks@\@xp{\shortauthors}\@temptokena\@xp{\shorttitle}%
  \toks4{\def\\{ \ignorespaces}}
  \edef\@tempa{%
    \@nx\markboth{\the\toks4
      \@nx\MakeUppercase{\the\toks@}}{\the\@temptokena}}%
  \@tempa
  \endgroup
  \c@footnote\z@
  \@cleartopmattertags
}
\makeatother

\title{On a conjecture of George Beck}

\author[S. Chern]{Shane Chern}
\address{Department of Mathematics, The Pennsylvania State University, University Park, PA 16802, USA}
\email{shanechern@psu.edu; chenxiaohang92@gmail.com}

\date{}

\begin{document}

{\footnotesize\noindent \textit{Int. J. Number Theory} \textbf{14} (2018), no.~3, 647--651. \MR{3786639}.\\
\doi{10.1142/S1793042118500392}}

\bigskip \bigskip

\maketitle

\begin{abstract}

In this paper, we prove a conjecture proposed by George Beck, which involves gap-free partitions and partitions with distinct parts.

\Keywords{Gap-free partition, distinct partition.}

\MSC{Primary 11P84; Secondary 05A17.}

\end{abstract}

\section{Introduction}

A \textit{partition} of a natural number $n$ is a nonincreasing sequence $\pi=(\pi_1,\pi_2,\ldots,\pi_\ell)$ of positive integers whose sum equals $n$, i.e. $\pi_1\ge \pi_2\ge\cdots\ge\pi_\ell$ and $\pi_1+\pi_2+\cdots+\pi_\ell=n$. Furthermore, a partition is called \textit{gap-free} if we restrict that the difference between each consecutive parts is at most $1$. For example, there are five gap-free partitions of $5$: $5$, $3+2$, $2+2+1$, $2+1+1+1$, and $1+1+1+1+1$. Let $a(n)$ denote the number of gap-free partitions of $n$. This sequence is listed as A034296 in the  On-Line Encyclopedia of Integer Sequence \cite{OEIS}. To determine the generating function of $a(n)$, we only need the following trivial observation (cf. \cite{And2017}):
\begin{quote}
The conjugates of gap-free partitions are partitions where only the largest part may repeat.
\end{quote}
Hence we have
\begin{equation}\label{eq:gen.a}
\sum_{n\ge 1}a(n)q^n=\sum_{m\ge 1}\frac{q^m}{1-q^m}(-q)_{m-1},
\end{equation}
where and in the sequel, we employ the standard $q$-series notations
$$(a)_n=(a;q)_n:=\prod_{k=0}^{n-1} (1-a q^k),$$
and
$$(a)_\infty=(a;q)_\infty:=\prod_{k\ge 0} (1-a q^k).$$

Recently, George Beck proposed an interesting conjecture in \cite[A034296]{OEIS}:

\begin{conjecture}\label{conj:beck}
$a(n)$ is also the sum of the smallest parts in the distinct partitions (i.e. partitions with distinct parts) of $n$ with an odd number of parts.
\end{conjecture}

Let $b(n)$ denote the sum of the smallest parts in the distinct partitions of $n$ with an odd number of parts. In this paper, we aim to prove Beck's conjecture.

\begin{theorem}\label{th:1}
For all $n\ge 1$, $a(n)=b(n)$.
\end{theorem}

\section{Proof of Theorem \ref{th:1}}

Let $\mathcal{D}$ be the set of distinct partitions.  For any $\pi\in\mathcal{D}$, we denote by $|\pi|$ the sum of the parts of $\pi$, by $\sigma(\pi)$ the smallest part of $\pi$, and by $\sharp(\pi)$ the number of parts of $\pi$. Then
$$\sum_{n\ge 1}b(n)q^n=\frac{1}{2}\left(\sum_{\pi\in \mathcal{D}}\sigma(\pi)q^{|\pi|}-\sum_{\pi\in \mathcal{D}}(-1)^{\sharp(\pi)}\sigma(\pi)q^{|\pi|}\right).$$
One also readily sees that
\begin{align*}
\sum_{\pi\in \mathcal{D}}\sigma(\pi)q^{|\pi|}&=\sum_{m\ge 1}mq^m(1+q^{m+1})(1+q^{m+2})\cdots\\
&=\sum_{m\ge 1}mq^m(-q^{m+1})_\infty,
\end{align*}
and
\begin{align*}
\sum_{\pi\in \mathcal{D}}(-1)^{\sharp(\pi)}\sigma(\pi)q^{|\pi|}&=\sum_{m\ge 1}(-m)q^m(1-q^{m+1})(1-q^{m+2})\cdots\\
&=-\sum_{m\ge 1}mq^m(q^{m+1})_\infty.
\end{align*}
Hence
\begin{align}
\sum_{n\ge 1}b(n)q^n&=\frac{1}{2}\left(\sum_{m\ge 1}mq^m(q^{m+1})_\infty+\sum_{m\ge 1}mq^m(-q^{m+1})_\infty\right)\notag\\
&=\frac{1}{2}\left((q)_\infty\sum_{m\ge 1}\frac{mq^m}{(q)_m}+(-q)_\infty\sum_{m\ge 1}\frac{mq^m}{(-q)_m}\right).\label{eq:2.1}
\end{align}

We first notice that
$$\sum_{m\ge 1}\frac{mq^m}{(q)_m}=\left.\left[z\frac{\partial}{\partial z}\left(\sum_{m\ge 0}\frac{z^m}{(q)_m}\right)\right]\right|_{z=q}.$$
Also the Euler's first sum \cite[Eq. (17.5.4)]{And2010} tells that
$$\sum_{m\ge 0}\frac{z^m}{(q)_m}=\frac{1}{(z)_\infty}.$$
Hence
\begin{align*}
z\frac{\partial}{\partial z}\left(\sum_{m\ge 0}\frac{z^m}{(q)_m}\right)&=z\frac{\partial}{\partial z}\left(\frac{1}{(z)_\infty}\right)\\
&=\frac{z}{(z)_\infty}\frac{\partial}{\partial z}\left(\log\frac{1}{(z)_\infty}\right)\\
&=\frac{z}{(z)_\infty}\frac{\partial}{\partial z}\left(\sum_{m\ge 0}\log\frac{1}{1-zq^m}\right)\\
&=\frac{z}{(z)_\infty}\sum_{m\ge 0}\frac{q^m}{1-zq^m}.
\end{align*}
This implies that
\begin{equation}\label{eq:2.2}
(q)_\infty\sum_{m\ge 1}\frac{mq^m}{(q)_m}=(q)_\infty\frac{q}{(q)_\infty}\sum_{m\ge 0}\frac{q^m}{1-q^{m+1}}=\sum_{m\ge 1}\frac{q^m}{1-q^m}.
\end{equation}

On the other hand, we have
$$\sum_{m\ge 1}\frac{mq^m}{(-q)_m}=\left.\left[z\frac{\partial}{\partial z}\left(\sum_{m\ge 0}\frac{z^m}{(-q)_m}\right)\right]\right|_{z=q}.$$
Let
$${}_{r+1}\phi_s\left(\begin{matrix} a_0,a_1,a_2\ldots,a_r\\ b_1,b_2,\ldots,b_s \end{matrix}; q, z\right):=\sum_{n\ge 0}\frac{(a_0;q)_n(a_1;q)_n\cdots(a_r;q)_n}{(q;q)_n(b_1;q)_n\cdots (b_s;q)_n}\left((-1)^n q^{\binom{n}{2}}\right)^{s-r}z^n.$$
We now need the Heine's first transformation \cite[Eq. (17.6.6)]{And2010}
$${}_{2}\phi_{1}\left(\begin{matrix} a,b\\ c \end{matrix}; q, z\right)=\frac{(b)_\infty(az)_\infty}{(c)_\infty(z)_\infty}{}_{2}\phi_{1}\left(\begin{matrix} c/b,z\\ az \end{matrix}; q, b\right).$$
Then
\begin{align*}
\sum_{m\ge 0}\frac{z^m}{(-q)_m}&={}_{2}\phi_{1}\left(\begin{matrix} 0,q\\ -q \end{matrix}; q, z\right)=\frac{(q)_\infty}{(-q)_\infty}\frac{1}{(z)_\infty}{}_{2}\phi_{1}\left(\begin{matrix} -1,z\\ 0 \end{matrix}; q, q\right)\\
&=\frac{(q)_\infty}{(-q)_\infty}\frac{1}{(z)_\infty}\sum_{n\ge 0}\frac{(-1)_n(z)_n}{(q)_n}q^n\\
&=\frac{(q)_\infty}{(-q)_\infty}\sum_{n\ge 0}\frac{(-1)_n}{(q)_n(zq^n)_\infty}q^n.
\end{align*}
It follows that
\begin{align*}
z\frac{\partial}{\partial z}\left(\sum_{m\ge 0}\frac{z^m}{(-q)_m}\right)&=\frac{z(q)_\infty}{(-q)_\infty}\sum_{n\ge 0}\frac{(-1)_nq^n}{(q)_n}\frac{\partial}{\partial z}\left(\frac{1}{(zq^n)_\infty}\right)\\
&=\frac{z(q)_\infty}{(-q)_\infty}\sum_{n\ge 0}\frac{(-1)_nq^n}{(q)_n}\frac{1}{(zq^n)_\infty}\sum_{m\ge 0}\frac{q^{m+n}}{1-zq^{m+n}}.
\end{align*}
Hence
\begin{align}
(-q)_\infty\sum_{m\ge 1}\frac{mq^m}{(-q)_m}&=(-q)_\infty\frac{q(q)_\infty}{(-q)_\infty}\sum_{n\ge 0}\frac{(-1)_nq^n}{(q)_n}\frac{1}{(q^{n+1})_\infty}\sum_{m\ge 0}\frac{q^{m+n}}{1-q^{m+n+1}}\notag\\
&=\sum_{n\ge 0}(-1)_nq^n\sum_{m\ge 0}\frac{q^{m+n+1}}{1-q^{m+n+1}}\notag\\
&=\sum_{m\ge 1}\frac{q^m}{1-q^m}\sum_{n=0}^{m-1}(-1)_nq^n\notag\\
&=\sum_{m\ge 1}\frac{q^m}{1-q^m}\left(1+2\sum_{n=1}^{m-1}(-q)_{n-1}q^n\right)\notag\\
&=\sum_{m\ge 1}\frac{q^m}{1-q^m}\left(2(-q)_{m-1}-1\right).\label{eq:2.3}
\end{align}
To see the last identity, we consider distinct partitions with largest part $\le m-1$. Obviously, its generating function is $(-q)_{m-1}$. On the other hand, if we fix the largest part to be $n\ge 1$, then it follows that
$$(-q)_{m-1}=1+\sum_{n=1}^{m-1}(-q)_{n-1}q^n.$$

Together with Eqs. \eqref{eq:2.1}--\eqref{eq:2.3}, we conclude that
\begin{align}
\sum_{n\ge 1}b(n)q^n&=\frac{1}{2}\left(\sum_{m\ge 1}\frac{q^m}{1-q^m}+\sum_{m\ge 1}\frac{q^m}{1-q^m}\left(2(-q)_{m-1}-1\right)\right)\notag\\
&=\sum_{m\ge 1}\frac{q^m}{1-q^m}(-q)_{m-1},
\end{align}
which coincides with the generating function of $a(n)$. Hence for all $n\ge 1$, $a(n)=b(n)$.

\section{Conclusion}

In a recent paper \cite{And2017b}, Andrews proved two conjectures of the same flavor (that is, showing some (weighted) partition functions are identical), which are also proposed by George Beck. We also notice the power of the differentiation technique in both our proof of Theorem \ref{th:1} and Andrews' proofs in \cite{And2017b}. Shishuo Fu and Dazhao Tang \cite{FT2017} subsequently generalized Andrews' results and provided a combinatorial proof. It would be interesting to know whether there also exists a combinatorial proof for our Theorem \ref{th:1}.

\bibliographystyle{amsplain}

\end{document}